\documentclass[11pt,reqno]{amsart}
\usepackage{graphicx}
\usepackage{verbatim}
\usepackage{textcomp}
\usepackage{amssymb}
\usepackage{cite}
\usepackage{amsmath}
\usepackage{latexsym}
\usepackage{amscd}
\usepackage{amsthm}
\usepackage{mathrsfs}
\usepackage{xypic}
\usepackage{bm}
\usepackage{url}

\newtheorem{thm}{Theorem}[section]

\newtheorem{lem}[thm]{Lemma}
\newtheorem{prop}[thm]{Proposition}

\theoremstyle{definition}

\theoremstyle{remark}

\numberwithin{equation}{section}
\def\f{\frac}

\def\ra{\rightarrow}
\def\pt{\partial}

\begin{document}
\title{Sharp diameter estimates for compact manifold with boundary}
\author{Haizhong Li}
\address{Department of mathematical sciences, and Mathematical Sciences
Center, Tsinghua University, 100084, Beijing, P. R. China}
\email{hli@math.tsinghua.edu.cn}
\author{Yong Wei}
\address{Department of mathematical sciences, Tsinghua University, 100084, Beijing, P. R. China}
\email{wei-y09@mails.tsinghua.edu.cn}
\thanks{The research of the authors was supported by NSFC No. 11271214.}

\maketitle

\begin{abstract}
Let $(N,g)$ be an $n$-dimensional complete Riemannian manifold with nonempty boundary $\pt N$. Assume that the Ricci curvature of $N$ has a negative lower bound $Ric\geq -(n-1)c^2$ for some $c>0$, and the mean curvature of the boundary $\pt N$ satisfies $H\geq (n-1)c_0>(n-1)c$ for some $c_0>c>0$. Then a known result (see \cite{LN}) says that $\sup_{x\in N}d(x,\pt N)\leq \frac 1c\coth^{-1}\frac{c_0}c$. In this paper, we prove that if the boundary $\pt N$ is compact, then the equality holds if and only if $N$ is isometric to the geodesic ball of radius $\frac 1c\coth^{-1}\frac{c_0}c$ in an $n$-dimensional hyperbolic space $\mathbb{H}^n(-c^2)$ of constant sectional curvature $-c^2$. Moreover, we also prove an analogous result for manifold with nonempty boundary and with $m$-Bakry-\'{E}mery Ricci curvature bounded below by a negative constant.
\end{abstract}

\section{Introduction}

Let $(N,g)$ be an $n$-dimensional complete Riemannian manifold. The classical Bonnet and Myers' theorem says that if the Ricci curvature of $(N,g)$ has a positive lower bound $Ric\geq (n-1)c^2>0$, then the diameter of $N$ is at most $\pi/{{c}}$. Cheng \cite{cheng} proved that if the diameter is equal to $\pi/{{c}}$, then $N$ is isometric to the $n$-sphere of constant sectional curvature $c^2$.

Recently, M. M. Li\cite{Martin} considered $n$-dimensional complete Riemannian manifold $(N,g)$ with nonnegative Ricci curvature and with mean convex boundary $\pt N$. M. M. Li proved that if the mean curvature of the boundary $\pt N$ satisifes $H\geq (n-1)c_0>0$ for some constant $c_0>0$, then
 \begin{eqnarray}\label{dist1}
    \sup_{x\in N}d(x,\pt N)&\leq &\f 1{c_0},
\end{eqnarray}
where $d$ denotes the distance function on $N$. Moreover, if $\pt N$ is compact, then $N$ is also compact and equality  holds in \eqref{dist1} if and only if $N$ is isometric to an $n$-dimensional Euclidean ball of radius $1/{c_0}$. Here the mean curvature $H$ of $\pt N$ is defined as the trace of the second fundamental form of $\pt N$ in $N$, that is, $H=\sum_{i=1}^{n-1}\langle\bar{\nabla}_{e_i}\nu,e_i\rangle$ for any orthonormal basis $e_1,\cdots,e_{n-1}$ of tangent bundle $T\pt N$, with respect to the outward unit normal $\nu$ of $\pt N$.

Note that by a similar argument as in the proof of the inequality \eqref{dist1}, if the Ricci curvature of $N$ has a negative lower bound $Ric\geq -(n-1)c^2$ for $c>0$, and the mean curvature of the  boundary $\pt N$ satisifes $H\geq (n-1)c_0>(n-1)c$ for some constant $c_0>c>0$, then one can prove that
\begin{eqnarray}\label{dist2}
    \sup_{x\in N}d(x,\pt N)&\leq&\frac 1{c}\coth^{-1}\frac{c_0}c.
\end{eqnarray}

In fact the distance bound \eqref{dist1} and \eqref{dist2} can be viewed as a Riemannian version of Hawking's singularity theorem (see e.g. \cite{HE}). The proof of the distance bound \eqref{dist1}, \eqref{dist2} is a standard argument by using the second variation formula of arc-length and can be found in Yanyan Li and Luc Nguyen's paper \cite[section 2]{LN}. In the first part of this paper, we study the equality case of \eqref{dist2}.  We have
\begin{thm}\label{main-thm1}
Let $(N^n,g)$ be an $n$-dimensional complete Riemannian manifold with nonempty boundary $\pt N$. Assume that the Ricci curvature of $N$ has a lower bound $Ric\geq -(n-1)c^2$ for some $c>0$, and the mean curvature of the  boundary $\pt N$ satisifes $H\geq (n-1)c_0>(n-1)c>0$ for some constant $c_0>c>0$. Then we know that the diameter estimate \eqref{dist2} holds in $N$.

If $\pt N$ is compact, then \eqref{dist2} implies that $N$ is also compact. We show if the equality holds in \eqref{dist2}, then $N$ is isometric to a geodesic ball of radius $\frac 1{c}\coth^{-1}\frac{c_0}c$ in an $n$-dimensional hyperbolic space $\mathbb{H}^n(-c^2)$ of constant sectional curvature $-c^2$.
\end{thm}

The proof of Theorem \ref{main-thm1} will be given in section 2. We first include the proof of the distance bound \eqref{dist2} for convenience of readers (see \cite{LN}). Then we consider the rigidity result when the equality occurs in \eqref{dist2}. By rescaling of metric, it suffices to consider the case $c=1$. The proof can be roughly divided into three parts. Firstly, by a Frankel type argument (see \cite{Law,Martin}), we show that under the curvature assumption of Theorem \ref{main-thm1}, the  boundary $\pt N$ is connected. Secondly, by using a similar argument as in M. M. Li's paper \cite{Martin}, we show that if the equality occurs in \eqref{dist2}, then $N$ is a geodesic ball of radius $\coth^{-1}c_0$ centered at some point $x_0$. Finally, by showing that the Laplacian comparison \eqref{lap-rho} assumes equality everywhere in $N$, we obtain that $N$ has constant sectional curvature $-1$ and is isometric to the hyperbolic ball. In the last step of the proof, a Heintze-Karcher's argument \cite{HK} will be used and is a key ingredient to the proof.

\vskip 2mm
Let $(N^n,g)$ be an $n$-dimensional complete smooth Riemannian manifold  and $f$ be a smooth function on $N$. We denote $\bar{\nabla},\bar{\Delta}$ and $\bar{\nabla}^2$ the gradient, Laplacian and Hessian operator on $N$ with respect to $g$, respectively. Given $m\in [n,\infty)$, the $m$-Bakry-\'{E}mery Ricci curvature of $(N,g)$ (see \cite{BM}) is defined by
\begin{align}\label{m-bakry}
    Ric_f^m=Ric+\bar{\nabla}^2f-\f 1{m-n}\bar{\nabla} f\otimes\bar{\nabla} f,\qquad (m> n).
\end{align}
When $m=\infty$, the last term of \eqref{m-bakry} is interpreted as 0 and this gives the Bakry-\'{E}mery Ricci curvature $ Ric_f=Ric+\bar{\nabla}^2f$. When $m=n$, this term only makes sense if $f$ is constant and in this case $Ric_f^m:=Ric$.

Recently, the study of  manifold with $m$-Bakry-\'{E}mery Ricci curvatures $Ric_f^m$ attracts many interests. Analogous to the Ricci curvature case, if one assumes that the $m$-Bakry-\'{E}mery Ricci curvature of $(N,g)$ satisfies $Ric_f^m\geq (m-1)c^2>0$, Qian \cite{Qian} proved that $diam(N)\leq \pi/{{c}}$ and then Ruan \cite{Ruan} proved that equality holds if and only if $N$ is isometric to the $n$-sphere of constant sectional curvature $c^2$. Recenlty, the authors \cite{Li-Wei} proved an analogous result of M. M. Li's theorem \cite{Martin} in the $m$-Bakry-\'{E}mery Ricci curvature case. Let $(N^n,g)$ be an $n$-dimensional complete Riemannian manifold with nonempty boundary $\pt N$ and $f$ be a smooth function on $N$. The $f$-mean curvature $H_f$ of $\pt N$ in $N$ is given by $H_f=H-\langle\bar{\nabla}f,\nu\rangle$, where $H$ is the mean curvature of $\pt N$ in N and $\nu$ is the outward unit normal of $\pt N$. When $f$ is constant, $H_f$ is just the mean curvature $H$. Assume that the $m$-Bakry-\'{E}mery Ricci curvature is nonnegative on $N$, and the $f$-mean curvature of the  boundary $\pt N$ satisifes $H_f\geq (m-1)c_0>0$ for some constant $c_0>0$. Then the authors \cite{Li-Wei} proved that
\begin{eqnarray}\label{dist3}
    \sup_{x\in N}d(x,\pt N)&\leq &\f 1{c_0}.
\end{eqnarray}
Moreover, if we assume that $\pt N$ is compact, then $N$ is also compact and equality  holds in \eqref{dist3} if and only if $N$ is isometric to an $n$-dimensional Euclidean ball of radius $1/{c_0}$.

In the second part of this paper, we consider the manifold with nonempty boundary and with $m$-Bakry-\'{E}mery Ricci curvature bounded below by a negative constant.

\begin{thm}\label{thm-2}
Let $(N^n,g)$ be an $n$-dimensional complete Riemannian manifold with nonempty boundary and $f$ be a smooth function on $N$. If the $m$-Bakry-\'{E}mery Ricci curvature of $N$ has a negative lower bound, i.e., $Ric_f^m\geq -(m-1)c^2$ for some constant $c>0$ on $N$, and the $f$-mean curvature of the  boundary $\pt N$ satisifes $H_f\geq (m-1)c_0>(m-1)c$ for some constant $c_0>c>0$, then we have the distance bound:
\begin{eqnarray}\label{dist4}
    \sup_{x\in N}d(x,\pt N)&\leq &\frac 1{c}\coth^{-1}\frac{c_0}c.
\end{eqnarray}
If $\pt N$ is compact, then \eqref{dist4} shows that $N$ is compact. Moreover, if the equality holds in \eqref{dist4}, we have $m=n$, and $N$ is isometric to a geodesic ball of radius $\frac 1{c}\coth^{-1}\frac{c_0}c$ in an $n$-dimensional hyperbolic space $\mathbb{H}^n(-c^2)$ of constant sectional curvature $-c^2$.
\end{thm}

The proof of Theorem \ref{thm-2} is similar with the proof of Theorem \ref{main-thm1} with some adjustment. Without loss of generality, we assume that $c=1$. When the equality occurs in \eqref{dist4}, arguing as the proof of Theorem \ref{main-thm1} we show that $N$ is equal to a geodesic ball of radius $\coth^{-1}c_0$ centered at some point $x_0$ and, the $f$-Laplacian comparison \eqref{f-Lap} assumes equality everywhere in $N$. There the generalized Heintze-Karcher theorem due to V. Bayle \cite{Bay} (see also \cite{Milman,Mor}) plays an important role. Finally, by using the Reilly formula for Bakry-\'{E}mery Ricci curvature (see \cite{MD,Li-Wei}), we show that $m=n$ and reduce to the case of Theorem \ref{main-thm1} and complete the proof of Theorem \ref{thm-2}.

\section{Proof of Theorem \ref{main-thm1}}

Firstly, for convenience of readers, we include the proof of \eqref{dist2} here (see \cite{LN}). For any point $x\in N$, since $N$ is complete, there exists a geodesic $\gamma:[0,d]\ra N$ parametrized by arc length with $\gamma(0)=x$, $\gamma(d)\in \pt N$ and $d=d(x,\pt N)$.  Choose  an orthonormal basis $e_1,\cdots e_{n-1}$ for $T_{\gamma(d)}\pt N$ and let $e_i(s)$ be the parallel transport of $e_i$ along $\gamma$. Let $V_i(s)=\varphi(s)e_i(s)$ with $\varphi(0)=0$ and $\varphi(d)=1$.  From the first variation formula, we have that for each $1\leq i\leq n-1$
\begin{align*}
    0=\delta\gamma(V_i)=&\langle \gamma'(d),V_i(d)\rangle-\langle \gamma'(0),V_i(0)\rangle-\int_0^d\langle{\gamma}''(s),V_i(s)\rangle ds\\
    =&\langle \gamma'(d),e_i(d)\rangle,
\end{align*}
which implies that $\gamma'(d)$ is orthogonal to $\pt N$ at $\gamma(d)$. The second variation formula gives that
\begin{align*}
    0\leq \sum_{i=1}^{n-1}\delta^2\gamma(V_i,V_i)=&\int_0^d\left((n-1)\varphi'(s)^2-\varphi(s)^2Ric(\gamma'(s),\gamma'(s))\right)ds\\
    &\quad +\langle \bar{\nabla}_{V_i(d)}V_i(d),\gamma'(d)\rangle-\langle \bar{\nabla}_{V_i(0)}V_i(0),\gamma'(0)\rangle\\
    =&\int_0^d\left((n-1)\varphi'(s)^2-\varphi(s)^2Ric(\gamma'(s),\gamma'(s))\right)ds-H(\gamma(d)).
\end{align*}
Since $Ric\geq -(n-1)c^2$ in $N$ and $H\geq (n-1)c_0>(n-1)c>0$ on $\pt N$, the above inequality implies
\begin{eqnarray}\label{eq1-0}
  0 &\leq& \int_0^d(\varphi'(s)^2+c^2\varphi^2(s))ds-c_0.
\end{eqnarray}
Choose
\begin{equation*}
    \varphi(s)=\frac{\sinh(cs)}{\sinh(cd)},\quad 0\leq s\leq d,
\end{equation*}
which satisfies $\varphi(0)=0$ and $\varphi(d)=1$. By substituting the above choosen $\varphi(s)$ into \eqref{eq1-0}, we have
\begin{eqnarray}
    c_0&\leq& c\coth(cd).
\end{eqnarray}
Therefore we have $d\leq \frac 1c\coth^{-1}\frac{c_0}c$ and this is the distance bound \eqref{dist2}.

The next lemma says that under the assumption of Theorem \ref{main-thm1}, the boundary $\pt N$ is connected.
\begin{lem}\label{lem2-1}
Let $(N^n,g)$ be an $n$-dimensional complete Riemannian manifold with nonempty boundary $\pt N$. Assume that the Ricci curvature of $N$ has a lower bound $Ric\geq -(n-1)c^2$, and the mean curvature of the  boundary $\pt N$ satisifes $H\geq (n-1)c_0>(n-1)c$ for some constant $c_0>c>0$. Then the boundary $\pt N$ is connected.
\end{lem}
\proof
We use the similar argument as in \cite{Law,Martin}. Suppose $\pt N$ is not connected, let $\Sigma$ be one of its components. Then $\Sigma$ and $\pt N\setminus\Sigma$ have a positive distance apart, i.e., $d(\Sigma,\pt N\setminus\Sigma)=l>0$. Since $\Sigma$ and $\pt N\setminus\Sigma$ are compact, there exists a minimizing geodesic $\gamma:[0,l]\rightarrow N$ parametrized by arc-length which realize the distance between $\Sigma$ and $\pt N\setminus\Sigma$. Note that $\gamma(0)\in\Sigma$, $\gamma(l)\in\pt N\setminus\Sigma$ and $\gamma(s)$ lies in the interior of $N$ for all $s\in(0,l)$. Moreover, $\gamma'(0)\bot T_{\gamma(0)}\pt N$ and $\gamma'(l)\bot T_{\gamma(l)}\pt N$, i.e.,$\gamma$ is a free boundary geodesic. Choose an orthonormal basis $e_1,\cdots,e_{n-1}$ for $T_{\gamma(0)}\pt N$ and let $e_i(s)$ be the parallel transport of $e_i$ along $\gamma$. Let $V_i(s)=\varphi(s)e_i(s)$ with $\varphi(0)=\varphi(l)=1$. Then the second variation formula of arc length gives that
\begin{align}
   0\leq\sum_{i=1}^{n-1}\delta^2\gamma(V_i,V_i)=&\int_0^l\left((n-1)\varphi'(s)^2-\varphi(s)^2Ric(\gamma'(s),\gamma'(s))\right)ds\nonumber\\
    &\qquad+\langle \bar{\nabla}_{V_i(l)}V_i(l),\gamma'(l)\rangle-\langle\bar{\nabla}_{V_i(0)}V_i(0),\gamma'(0)\rangle\nonumber\\
    =&\int_0^l\left((n-1)\varphi'(s)^2-\varphi(s)^2Ric(\gamma'(s),\gamma'(s))\right)ds\nonumber\\
    &\qquad-H(\gamma(l))-H(\gamma(0)).\label{eq1-1}
\end{align}
Now we choose
\begin{equation*}
  \varphi(s)=\frac{\cosh c(s-\frac l2)}{\cosh (cl/2)}, \quad 0\leq s\leq l.
\end{equation*}
Then $\varphi(0)=\varphi(l)=1$. Substituting the above chosen $\varphi(s)$ into \eqref{eq1-1},and using the assumption $Ric\geq -(n-1)c^2$ in $N$ and $H\geq (n-1)c_0>(n-1)c>0$ on $\pt N$,  we obtain
\begin{eqnarray}\label{eq1-2}
  0 &\leq& 2(n-1)c\tanh(\frac{cl}2)-2(n-1)c_0.
\end{eqnarray}
Since $l$ is the distance of $\Sigma$ and $\pt N\setminus\Sigma$ and $\gamma(t)$ realizes this distance, we have $l\leq 2\sup_{x\in N}d(x,\pt N)\leq \frac 2c\coth^{-1}\frac{c_0}c$ by \eqref{dist2}. Therefore the right-hand side of \eqref{eq1-2} is bounded from above by
\begin{eqnarray*}
  && 2(n-1)c\tanh(\coth^{-1}\frac{c_0}c)-2(n-1)c_0 \\
  &=&\frac{2(n-1)}{c_0}(c^2-c_0^2)<0,
\end{eqnarray*}
which is a contradiction. We conclude that the boundary $\pt N$ is connected.
\endproof

Now we continue the proof of Theorem \ref{main-thm1}. Assume that $\pt N$ is compact, then \eqref{dist2} implies $N$ is also compact. Suppose equality holds in \eqref{dist2}. By rescaling the metric of $N$, we can  assume that $c=1$. Since $N$ is compact, there exists one point $x_0\in N$ such that
 \begin{equation}\label{eq1-3}
   d(x_0,\pt N)=\coth^{-1}c_0.
 \end{equation}
 We denote $\rho_0=\coth^{-1}c_0$ for simplicity.  We first show that
 \begin{lem}\label{lem2-2}
 Under the assumption of Theorem \ref{main-thm1}, if the equality holds in \eqref{dist2}, then $N$ is equal to the geodesic ball of radius $\rho_0$ centered at $x_0$.
 \end{lem}
\proof
 It is clear that the geodesic ball $B_{\rho_0}(x_0)$ of radius $\rho_0$ centered at $x_0$ is contained in $N$. We show that $N$ is just equal to the geodesic ball $B_{\rho_0}(x_0)$. Let $\rho=d(x_0,\cdot)$ be the distance function from $x_0$. Since the Ricci curvature of $N$ satisfies $Ric\geq -(n-1)$, the Laplacian of $\rho$ satisfies (see \cite{Schoen-Yau,peter})
\begin{eqnarray}\label{lap-rho}
    \bar{\Delta}\rho&\leq&(n-1)\coth\rho,
\end{eqnarray}
in the sense of distribution. Let $\Sigma=\{q\in\pt N:\rho(q)=\rho_0\}$, which is clearly a closed set in $\pt N$ by the continuity of $\rho$. Since $\pt N$ is connected, to show $\Sigma=\pt N$, it suffices to show that $\Sigma$ is also open in $\pt N$, that is for any $q\in \Sigma$, there is an open neighborhood $U$ of $q$ in $\pt N$ such that $\rho\equiv \rho_0$ on $U$. If $q$ is not a conjugate point to $x_0$ in $N$, then the geodesic sphere $\pt B_{\rho_0}(x_0)$ is a smooth hypersurface near $q$ in $N$. Note that $\bar{\Delta}\rho$ is the mean curvature $H$ of the geodesic sphere and $\rho=\rho_0$ on the geodesic sphere, by the Laplacian comparison inequality \eqref{lap-rho}, we have that the mean curvature $H$ of the geodesic sphere is at most $(n-1)c_0$. However, by the assumption of Theorem \ref{main-thm1}, the mean curvature of $\pt N$ is at least $(n-1)c_0$. Then from the maximum principle (see \cite{Esch}), we have that $\pt N$ and $\pt B_{\rho_0}(x_0)$ coincides in a neighborhood of $q$.  This implies that $\Sigma$ is open near any $q$ which is not a conjugate point. A similar process in \cite{Martin} (see also Calabi \cite{Cal}) makes us to work through the argument to conclude that $\rho$ is constant near $q$ in $\pt N$, when $q$ is a conjugate point of $x_0$. This proves that $\pt N$ is just the geodesic sphere $\pt B_{\rho_0}(x_0)$ and $N$ is the geodesic ball $B_{\rho_0}(x_0)$ of radius $\rho_0$ centered at $x_0$.
\endproof

Since any $q\in \pt N$ can be joined by a minimizing geodesic $\gamma$ parameterized by arc-length from $x_0$ to $q$, and $\gamma$ is orthogonal to $\pt N=\pt B_{\rho_0}(x_0)$ at $q$, the geodesic $\gamma$ is uniquely determined by $q$ and $q$ is not in the cut locus of $x_0$. So that $\rho=d(x_0,\cdot)$ is smooth up to the boundary $\pt N$. From the proof of lemma \ref{lem2-2}, the mean curvature $H$ of the boundary $\pt N$ satisfies $H=(n-1)\coth\rho_0$. Since $\pt N$ is the geodesic sphere $\pt B_{\rho_0}(x_0)$, the Laplacian of $\rho$ is equal to the mean curvature $H$ on the boundary $\pt N$. Therefore \begin{eqnarray}\label{lap-rho1}
    \bar{\Delta}\rho&=&(n-1)\coth\rho,
\end{eqnarray}
holds on the boundary $\pt N$. We next show that \eqref{lap-rho1} holds everywhere in $N$.

On the one hand, since $Ric\geq -(n-1)$ on $N$, the Laplacian comparison \eqref{lap-rho} for $\rho$ holds in the classical sense. Note that $|\bar{\nabla}\rho|=1$ in $N$, then we have
\begin{eqnarray*}
  \bar{\Delta}\cosh\rho &=& \bar{\Delta}\rho\sinh\rho+\cosh\rho|\bar{\nabla}\rho|^2\leq n\cosh\rho.
\end{eqnarray*}
Integrating the above inequality over $N$ and by divergence theorem, we have
\begin{eqnarray*}
  \int_{\pt N}\sinh\rho\langle\bar{\nabla}\rho,\nu\rangle d\mu&\leq & n\int_N\cosh\rho dv,
\end{eqnarray*}
where $\nu$ is the outward unit normal of $\pt N$ in $N$ and $d\mu$, $dv$ are volume elements on $\pt N$ and $N$ respectively. Note that $\rho=\rho_0$ and $\langle\bar{\nabla}\rho,\nu\rangle=1$ on $\pt N$, the above inequality gives that
\begin{eqnarray}\label{eq1-6}
  \int_{\pt N}\sinh\rho_0 d\mu&\leq &n\int_N\cosh\rho(x)dv.
\end{eqnarray}

On the other hand, we prove the reversed inequality also holds in \eqref{eq1-6}.
\begin{lem}
We have
\begin{eqnarray}\label{eq1-7}
  \int_{\pt N}\sinh\rho_0 d\mu&\geq& n\int_N\cosh \rho(x)dv.
\end{eqnarray}
\end{lem}
\proof
The proof is inspired by Heintze-Karcher's paper \cite{HK}. Note that any $q\in \pt N$ can be joined by a minimizing geodesic $\gamma$ parameterized by arc-length from $x_0$ to $q$, and $\gamma$ is orthogonal to $\pt N=\pt B_{\rho_0}(x_0)$ at $q$, the geodesic $\gamma$ is uniquely determined by $q$ and $q$ is not in the cut locus of $x_0$. The exponential map of the normal bundle $T^{\bot}\pt N$ of $\pt N$ in $N$ is surjective.  For any $q\in\pt N$, the curve $\gamma(t)=\exp_q(-t\nu)$ $(0\leq t\leq\rho_0)$ is the geodesic connecting $q$ and $x_0$, i.e., $\gamma(\rho_0)=x_0$. We have
\begin{eqnarray*}
\int_N\cosh \rho(x) dv&\leq& \int_{\pt N}\int_0^{\rho_0}\cosh \rho(\exp_q(-t\nu(q))) |\det(d \exp_q)_{t,\nu}|dtd\mu_q,
\end{eqnarray*}
Since $Ric\geq -(n-1)$ in $N$, Corollary 3.3.2 of \cite{HK} gives that
\begin{eqnarray}
|\det(d \exp_q)_{t\nu}| &\leq&(\cosh t-\frac{H(q)}{n-1}\sinh t)^{n-1}.
\end{eqnarray}
Note that $\rho(\exp_q(-t\nu(q)))=\rho_0-t$, and
\begin{eqnarray*}
    \cosh (\rho_0-t)&=&\cosh\rho_0\cosh t-\sinh\rho_0\sinh t\\
    &=&\sinh\rho_0(\coth\rho_0\cosh t-\sinh t).
\end{eqnarray*}
From the proof of Lemma \ref{lem2-2}, we have $H(q)=(n-1)\coth\rho_0$ on $\pt N$. So we have
\begin{eqnarray*}
\int_N\cosh \rho(x) dv&\leq& -\int_{\pt N}\int_0^{\rho_0}\frac 1n\sinh\rho_0 \frac d{dt}(\cosh t-\coth\rho_0\sinh t)^{n}dtd\mu_q\\
&=&\frac 1n\int_{\pt N}\sinh\rho_0d\mu
\end{eqnarray*}
which gives the inequality \eqref{eq1-7}.
\endproof

Combining \eqref{eq1-6} and \eqref{eq1-7}, we have that \eqref{lap-rho1} holds everywhere in $N$. Thus we conclude that $N$ has constant sectional curvature $-1$ and is isometric to the hyperbolic ball (see \cite{peter}). The proof of Theorem \ref{main-thm1} is completed.

\section{Proof of Theorem \ref{thm-2}}

Firstly, the proof of diameter estimate \eqref{dist4} is similar with the proof of \eqref{dist2}, which is also by using the second variation formula of arc-length. For any point $x\in N$, since $N$ is complete, there exists a geodesic $\gamma:[0,d]\ra N$ parametrized by arc-length with $\gamma(0)=x$, $\gamma(d)\in \pt N$ and $d=d(x,\pt N)$.  Choose  an orthonormal basis $e_1,\cdots e_{n-1}$ for $T_{\gamma(d)}\pt N$ and let $e_i(s)$ be the parallel transport of $e_i$ along $\gamma$. Let $V_i(s)=\varphi(s)e_i(s)$ with $\varphi(0)=0$ and $\varphi(d)=1$. The first variation formula implies that $\gamma'(d)$ is orthogonal to $\pt N$ at $\gamma(d)$. The second variation formula gives that
\begin{align*}
    0\leq \sum_{i=1}^{n-1}\delta^2\gamma(V_i,V_i)=&\int_0^d\left((n-1)\varphi'(s)^2-\varphi(s)^2Ric(\gamma'(s),\gamma'(s))\right)ds-H(\gamma(d)).
\end{align*}
By the definition of $m$-Bakry-\'{E}mery Ricci curvature, we have
\begin{align}
    0\leq&\int_0^d\left((n-1)\varphi'(s)^2-\varphi(s)^2Ric_f^m(\gamma'(s),\gamma'(s))\right)ds-H(\gamma(d))\nonumber\\
    &\qquad+\int_0^d\varphi(s)^2\left(\bar{\nabla}^2f(\gamma'(s),\gamma'(s))-\f 1{m-n}\langle\bar{\nabla}f(\gamma(s)),\gamma'(s)\rangle^2\right)ds\nonumber\\
    =&\int_0^d\left((n-1)\varphi'(s)^2-\varphi(s)^2Ric_f^m(\gamma'(s),\gamma'(s))\right)ds-H(\gamma(d))\nonumber\\
    &\qquad+\int_0^d\varphi(s)^2\left(\f{d^2}{ds^2}f(\gamma(s))-\f 1{m-n} (\f{d}{ds}f(\gamma(s)))^2\right)ds\label{eq3-2}
\end{align}
where we used the facts
\begin{align*}
    \f{d}{ds}f(\gamma(s))=\langle\bar{\nabla}f(\gamma(s)),\gamma'(s)\rangle
\end{align*}
and
\begin{align*}
    \f{d^2}{ds^2}f(\gamma(s))=&\bar{\nabla}^2f(\gamma'(s),\gamma'(s)).
\end{align*}
By integration by parts, we deduce from \eqref{eq3-2} that
\begin{align}
    0\leq&\int_0^d\biggl((n-1)\varphi'(s)^2-\varphi(s)^2Ric_f^m(\gamma'(s),\gamma'(s))-2\varphi(s)\varphi'(s)\f{d}{ds}f(\gamma(s))\nonumber\\
    &\qquad-\f 1{m-n}\varphi(s)^2(\f{d}{ds}f(\gamma(s)))^2\biggr)ds+\varphi(d)^2\langle\bar{\nabla}f(\gamma(d)),\gamma'(d)\rangle\nonumber\\
    &\qquad -\varphi(0)^2\langle\bar{\nabla}f(\gamma(0)),\gamma'(0)\rangle-H(\gamma(d)).\label{eq3-3}
\end{align}
Note that $\varphi(0)=0, \varphi(d)=1$ and $\gamma'(d)$ is equal to the outward unit normal vector $\nu$ at $\gamma(d)\in\pt N$. The $f$-mean curvature $H_f$ at $\gamma(d)$ is
\begin{eqnarray*}
H_f(\gamma(d))&=& H(\gamma(d))-\langle\bar{\nabla}f(\gamma(d)),\nu(\gamma(d))\rangle.
\end{eqnarray*}
Moreover, the Cauchy-Schwartz inequality implies
\begin{eqnarray*}
  -2\varphi(s)\varphi'(s)\f{d}{ds}f(\gamma(s))&\leq&(m-n)\varphi'(s)^2+\f 1{m-n}\varphi(s)^2(\f{d}{ds}f(\gamma(s)))^2.
\end{eqnarray*}
Thus from \eqref{eq3-3}, we have
\begin{align}\label{eq2-1}
    0\leq&\int_0^d\biggl((m-1)\varphi'(s)^2-\varphi(s)^2Ric_f^m(\gamma'(s),\gamma'(s))\biggr)ds-H_f(\gamma(d)).
\end{align}
Choose
\begin{equation*}
    \varphi(s)=\frac{\sinh(cs)}{\sinh(cd)},\quad 0\leq s\leq d,
\end{equation*}
which satisfies $\varphi(0)=0$ and $\varphi(d)=1$. Since $Ric_f^m\geq -(m-1)c^2$ in $N$ and $H_f\geq (m-1)c_0>(m-1)c$ on $\pt N$, by substituting the above choosen $\varphi(s)$ into \eqref{eq2-1}, we have
\begin{eqnarray}
    c_0&\leq& c\coth(cd).
\end{eqnarray}
Therefore we have $d\leq \frac 1c\coth^{-1}\frac{c_0}c$ and this is the distance bound \eqref{dist4}.

\vskip 2mm
If the boundary $\pt N$ is compact, then \eqref{dist4} implies $N$ is also compact. Next we prove the rigidity result when the equality occurs in \eqref{dist4}. As in the proof of Theorem \ref{main-thm1}, the curvature assumption of Theorem \ref{thm-2} implies the boundary $\pt N$ is connected.
\begin{lem}\label{lem3-1}
Under the curvature assumption of Theorem \ref{thm-2}, the boundary $\pt N$ is connected.
\end{lem}
\proof
The proof is also by a Frankel type argument, see lemma \ref{lem2-1}. We include a proof here for exhibiting the adjustment. Suppose $\pt N$ is not connected, let $\Sigma$ be one of its components. Let $\gamma(s)$ ($0\leq s\leq l$) be the free boundary geodesic realizing the distance between $\Sigma$ and $\pt N\setminus\Sigma$. Choose an orthonormal basis $e_1,\cdots,e_{n-1}$ for $T_{\gamma(0)}\pt N$ and let $e_i(s)$ be the parallel transport of $e_i$ along $\gamma$. Let $V_i(s)=\varphi(s)e_i(s)$ with $\varphi(0)=\varphi(l)=1$. Then the second variation formula of arc-length gives that
\begin{align*}
   0\leq\sum_{i=1}^{n-1}\delta^2\gamma(V_i,V_i)=&\int_0^l\left((n-1)\varphi'(s)^2-\varphi(s)^2Ric(\gamma'(s),\gamma'(s))\right)ds\\
    &\qquad+\langle \bar{\nabla}_{V_i(l)}V_i(l),\gamma'(l)\rangle-\langle\bar{\nabla}_{V_i(0)}V_i(0),\gamma'(0)\rangle\\
    =&\int_0^l\left((n-1)\varphi'(s)^2-\varphi(s)^2Ric(\gamma'(s),\gamma'(s))\right)ds\\
    &\qquad-H(\gamma(l))-H(\gamma(0)).
\end{align*}
By the definition of $m$-Bakry-\'{E}mery Ricci curvature, and using the Cauchy-Schwartz inequality as in the proof of \eqref{dist4}, we have
\begin{eqnarray}
0 &\leq& \int_0^l\biggl((m-1)\varphi'(s)^2-\varphi(s)^2Ric_f^m(\gamma'(s),\gamma'(s))\biggr)ds\nonumber\\
 &&\qquad-H_f(\gamma(l))-H_f(\gamma(0)).\label{eq3-1}
\end{eqnarray}
Since $Ric_f^m\geq -(m-1)c^2$ in $N$ and $H_f\geq (m-1)c_0>(m-1)c>0$ on $\pt N$, we can argue as the proof of lemma \ref{lem2-1} to get a contradiction by choosing the function
\begin{equation*}
  \varphi(s)=\frac{\cosh c(s-\frac l2)}{\cosh (cl/2)}, \quad 0\leq s\leq l.
\end{equation*}
in \eqref{eq3-1}. Then we conclude that the boundary $\pt N$ is connected.
\endproof

Now assume that the equality occurs in \eqref{dist4}. Without loss of generality, we assume that $c=1$. By the compactness of $N$, there exists one point $x_0\in N$ such that
\begin{eqnarray}
d(x_0,\pt N)&=&\coth^{-1}c_0.
\end{eqnarray}
We also denote $\rho_0=\coth^{-1}c_0$ for simplicity.
\begin{lem}\label{lem3-2}
Under the assumption of Theorem \ref{thm-2}, if the equality holds in \eqref{dist4}, then $N$ is equal to the geodesic ball of radius $\rho_0$ centered at $x_0$.
\end{lem}
\proof
The proof is similar with the proof of lemma \ref{lem2-2}. The only difference is that we replace the Laplacian comparison \eqref{lap-rho} by the following $f$-Laplacian comparison. Since $Ric_f^m\geq -(m-1)$ in $N$, the $f$-Laplacian comparison of the distance function $\rho(x)=d(x_0,x)$ due to Qian \cite{Qian} says that
\begin{eqnarray}\label{f-Lap}
\bar{\Delta}_f\rho(x):=\bar{\Delta}\rho(x)-\bar{\nabla}f\cdot\bar{\nabla}\rho(x) &\leq& (m-1)\coth\rho(x).
\end{eqnarray}
holds in the sense of distribution.
\endproof

Next we show that the $f$-Laplacian comparison \eqref{f-Lap} assumes equality everywhere in $N$.  From lemma \ref{lem3-2}, $N$ is the geodesic ball of radius $\rho_0$ centered at $x_0$.  Any $q\in \pt N$ can be joined by a minimizing geodesic $\gamma$ parameterized by arc-length from $x_0$ to $q$, and $\gamma$ is orthogonal to $\pt N=\pt B_{\rho_0}(x_0)$ at $q$. The geodesic $\gamma$ is uniquely determined by $q$ and $q$ is not in the cut locus of $x_0$. Then the distance function $\rho(x)$ is smooth up to the boundary $\pt N$. The $f$-Laplacian comparison \eqref{f-Lap} implies
\begin{eqnarray*}
  \bar{\Delta}_f\cosh\rho(x) &=& \bar{\Delta}_f\rho(x)\sinh\rho(x)+\cosh\rho(x)|\bar{\nabla}\rho(x)|^2\leq m\cosh\rho(x).
\end{eqnarray*}
Integrating the above inequality over $N$ with respect to the weighted volume element $e^{-f}dv$ and by divergence theorem, we have
\begin{eqnarray*}
  \int_{\pt N}\sinh\rho\langle\bar{\nabla}\rho,\nu\rangle e^{-f}d\mu&\leq & m\int_N\cosh\rho(x) ~e^{-f}dv,
\end{eqnarray*}
where $\nu$ is the outward unit normal of $\pt N$ in $N$. Note that $\rho=\rho_0$ and $\langle\bar{\nabla}\rho,\nu\rangle=1$ on $\pt N$, the above inequality gives that
\begin{eqnarray}\label{eq2-6}
  \int_{\pt N}\sinh\rho_0 e^{-f}d\mu&\leq &m\int_N\cosh\rho(x)e^{-f}dv.
\end{eqnarray}

On the other hand, we prove the reversed inequality also holds in \eqref{eq2-6}.
\begin{lem}
We have
\begin{eqnarray}\label{eq2-7}
  \int_{\pt N}\sinh\rho_0  e^{-f}d\mu&\geq& m\int_N\cosh \rho(x)e^{-f}dv.
\end{eqnarray}
\end{lem}
\proof
To show \eqref{eq2-7}, we need the generalized Heintze-Karcher theorem due to V. Bayle \cite{Bay}. Note that any $q\in \pt N$ can be joined by a minimizing geodesic $\gamma$ parameterized by arc-length from $x_0$ to $q$, and $\gamma$ is orthogonal to $\pt N=\pt B_{\rho_0}(x_0)$ at $q$, the geodesic $\gamma$ is uniquely determined by $q$ and $q$ is not in the cut locus of $x_0$. The exponential map of the normal bundle $T^{\bot}\pt N$ of $\pt N$ in $N$ is surjective.  For any $q\in\pt N$, the curve $\gamma(t)=\exp_q(-t\nu)$ $(0\leq t\leq\rho_0)$ is the geodesic connecting $q$ and $x_0$, i.e., $\gamma(\rho_0)=x_0$. Since $Ric_f^m\geq -(m-1)$ in $N$,   the generalized Heintze-Karcher theorem in \cite{Bay} implies
\begin{eqnarray*}
&&\int_N\cosh \rho(x) e^{-f}dv\\
&\leq& \int_{\pt N}\int_0^{\rho_0}\cosh \rho(\exp_q(-t\nu(q)))(\cosh t-\frac{H_f(q)}{m-1}\sinh t)^{m-1}dte^{-f(q)}d\mu_q,
\end{eqnarray*}
Note that $\rho(\exp_q(-t\nu(q)))=\rho_0-t$ and
\begin{eqnarray*}
    \cosh (\rho_0-t)&=&\cosh\rho_0\cosh t-\sinh\rho_0\sinh t\\
    &=&\sinh\rho_0(\coth\rho_0\cosh t-\sinh t).
\end{eqnarray*}
From the proof of lemma \ref{lem3-2} we have $H_f(q)=(m-1)\coth\rho_0$ on $\pt N$. Therefore
\begin{eqnarray*}
&&\int_N\cosh \rho(x) e^{-f}dv\\
&\leq& -\int_{\pt N}\int_0^{\rho_0}\frac 1m\sinh\rho_0 \frac d{dt}(\cosh t-\coth\rho_0\sinh t)^{m}dte^{-f(q)}d\mu_q\\
&=&\frac 1m\int_{\pt N}\sinh\rho_0e^{-f}d\mu
\end{eqnarray*}
which gives the inequality \eqref{eq2-7}.
\endproof

Combining \eqref{eq2-6} and \eqref{eq2-7}, the $f$-Laplacian comparison inequality \eqref{f-Lap} assumes equality everywhere in $N$, i.e., we have that
\begin{eqnarray}\label{f-Lap-eq}
\bar{\Delta}_f\rho(x) &=& (m-1)\coth\rho(x)
\end{eqnarray}
holds in the classical sense everywhere in $N$. Finally, we show that $m=n$.

\begin{lem}
We have $m=n$.
\end{lem}
\proof
Recall that for any function $u\in C^3(N)$, Ma-Du \cite{MD} obtained the following Reilly formula for Bakry-\'{E}mery Ricci curvature $Ric_f$:
\begin{align}
    0=&\int_{N}(Ric_f(\bar{\nabla}u,\bar{\nabla}u)-|\bar{\Delta}_fu|^2+|\bar{\nabla}^2u|^2)e^{-f}dv\label{Reilly}\\
    &\quad +\int_{\pt N}\left((\Delta_fu+H_f\f{\pt u}{\pt \nu})\f{\pt u}{\pt \nu}-\langle \nabla u,\nabla\frac{\pt u}{\pt\nu}\rangle+h(\nabla u,\nabla u)\right)e^{-f}d\mu.\nonumber
\end{align}
Here, $\bar{\Delta}_f=\bar{\Delta}-\bar{\nabla}f\cdot\bar{\nabla},\bar{\nabla}$ and $\bar{\nabla}^2$ are the $f$-Laplacian, gradient and Hessian on $N$ respectively; $\Delta_f=\Delta-\nabla f\cdot \nabla $ and $\nabla$ are the $f$-Laplacian and gradient operators on $\pt N$; $\nu$ is the outward unit normal of $\pt N$; $H_f$ and $h$ are the $f$-mean curvature and second fundamental form of $\pt N$ in $N$ with respect to $\nu$ respectively.

Suppose on the contrary we have $m>n$, let $z=\f mn$ and by the basic algebraic inequality $(a+b)^2\geq \f{a^2}{z}-\f{b^2}{z-1}$ for $z>1$, we have (see \cite{XDLi,Li-Wei})
\begin{eqnarray}
    |\bar{\nabla}^2u|^2\geq \f 1n(\bar{\Delta} u)^2&=&\f 1n(\bar{\Delta}_fu+\bar{\nabla}f\cdot\bar{\nabla}u)^2\nonumber\\
    &\geq&\f 1n\left(\f nm(\bar{\Delta}_fu)^2-\f n{m-n}(\bar{\nabla}f\cdot\bar{\nabla}u)^2\right)\nonumber\\
    &=&\f 1m(\bar{\Delta}_fu)^2-\f 1{m-n}(\bar{\nabla}f\cdot\bar{\nabla}u)^2.\label{eq2-alg}
\end{eqnarray}
Substituting this into \eqref{Reilly} and using the definition \eqref{m-bakry} of $m$-Bakry-\'{E}mery Ricci curvature, we get
\begin{align}
    0\geq&\int_{N}(Ric_f^m(\bar{\nabla}u,\bar{\nabla}u)-\f {m-1}m|\bar{\Delta}_fu|^2)e^{-f}dv\label{Reilly-2}\\
    &\quad +\int_{\pt N}\left((\Delta_fu+H_f\f{\pt u}{\pt \nu})\f{\pt u}{\pt \nu}-\langle \nabla u,\nabla\frac{\pt u}{\pt\nu}\rangle+h(\nabla u,\nabla u)\right)e^{-f}d\mu.\nonumber
\end{align}

Now since  \eqref{f-Lap-eq} holds everywhere in $N$, we have that
\begin{eqnarray}\label{eq2-8}
  \bar{\Delta}_f\cosh\rho(x) &=&m\cosh\rho(x)
\end{eqnarray}
also holds everywhere in $N$. Note that $\pt N$ is a geodesic ball of radius $\rho_0$ centered at $x_0$, $Ric_f^m\geq -(m-1)$ in $N$ and $H_f=(m-1)c_0=(m-1)\coth\rho_0$ on $\pt N$. Substituting $u(x)=\cosh\rho(x)$ into the Reilly inequality \eqref{Reilly-2}, using \eqref{eq2-8} and integrating by part, we obtain
\begin{eqnarray*}
0&\geq&\int_N\big(-(m-1)|\bar{\nabla}u|^2-\f {m-1}m|\bar{\Delta}_fu|^2\big)e^{-f}dv+\int_{\pt N}H_f(\f{\pt u}{\pt \nu})^2e^{-f}d\mu\\
&=&\frac{m-1}m\int_N\bar{\Delta}_fu\left(mu-\bar{\Delta}_fu\right)e^{-f}dv-\int_{\pt N}(m-1)\frac{\pt u}{\pt \nu}ue^{-f}d\mu\\
&& \qquad +(m-1)\int_{\pt N}\coth\rho_0(\frac{\pt u}{\pt \nu})^2e^{-f}d\mu\\
&=&-\int_{\pt N}(m-1)\sinh\rho_0\cosh\rho_0e^{-f}d\mu+(m-1)\int_{\pt N}\coth\rho_0(\sinh\rho_0)^2e^{-f}d\mu\\
&=&0,
\end{eqnarray*}
where we used the facts $\rho=\rho_0$ and $\frac{\pt \rho}{\pt \nu}=1$ on $\pt N$. Therefore, the algebraic inequality \eqref{eq2-alg} assumes equality everywhere for $u(x)=\cosh\rho(x)$. Thus we have
\begin{eqnarray*}
  0&=&\bar{\Delta}_f\cosh\rho(x)+\frac m{m-n}\bar{\nabla}f\cdot\bar{\nabla}\cosh\rho(x)\\
  &=&\bar{\Delta}\cosh\rho(x)+\frac n{m-n}\bar{\nabla}f\cdot\bar{\nabla}\cosh\rho(x)
\end{eqnarray*}
holds everywhere in $N$. Let $\omega(x)=\cosh\rho(x)-\cosh\rho_0$. Then
\begin{eqnarray}
0&=&\bar{\Delta}\omega(x)+\frac n{m-n}\bar{\nabla}f\cdot\bar{\nabla}\omega(x)\label{eq2-9}
\end{eqnarray}
in $N$ and $\omega(x)=0$ on $\pt N$. Multiplying \eqref{eq2-9} with $\omega(x)$ and integrating over $N$ with respect to $e^{\frac n{m-n}f}dv$, we get
\begin{eqnarray*}
0 &=& \int_N\omega\bigl(\bar{\Delta}\omega(x)+\frac n{m-n}\bar{\nabla}f\cdot\bar{\nabla}\omega(x)\bigr)e^{\frac n{m-n}f}dv\\
&=&-\int_N|\bar{\nabla}\omega|^2e^{\frac n{m-n}f}dv+\int_{\pt N}\omega\frac{\pt\omega}{\pt\nu}e^{\frac n{m-n}f}d\mu\\
&=&-\int_N|\bar{\nabla}\omega|^2e^{\frac n{m-n}f}dv
\end{eqnarray*}
where the third equality is due to the fact $\omega(x)=0$ on $\pt N$. Therefore we have that $\omega(x)=\cosh\rho(x)-\cosh\rho_0$ is constant in $N$, which is a contradiction. Thus we conclude that $m=n$.
\endproof

Once we have $m=n$, the last statement of Theorem \ref{thm-2} follows from Theorem \ref{main-thm1}, and we complete the proof of Theorem \ref{thm-2}.

\appendix
\section{Manifold with $Ric_f$ bounded below}

In this appendix, we give a result on the diameter estimate for manifold $(N^n,g)$ with Bakry-\'{E}mery Ricci curvature $Ric_f$ bounded below. By assuming that the function $f$ is bounded, i.e., $|f|\leq k$, and $Ric_f\geq (n-1)c^2>0$ in $N$, Wei-Wylie \cite{Wei-Wylie} proved that the diameter of $N$ satisfies $diam(N)\leq (\pi+\frac{4k}{n-1})/{{c}}$. See a different upper bound $diam(N)\leq \sqrt{1+\frac{2\sqrt{2}k}{n-1}}\pi/{{c}}$ obtained by Limoncu \cite{Limoncu}. The following proposition deals with the manifold with $Ric_f\geq -(n-1)c^2$ for some $c\geq 0$ and with nonempty boundary. By assuming $|f|\leq k$ for some constant $k$, we have

\begin{prop}\label{prop-5}
Let $(N^n,g)$ be an $n$-dimensional complete Riemannian manifold with nonempty boundary and $f$ be a smooth bounded function ($|f|\leq k$) on $N$. Assume that the Bakry-\'{E}mery Ricci curvature $Ric_f\geq -(n-1)c^2$ for some $c\geq 0$ on $N$, and the $f$-mean curvature of the  boundary $\pt N$ satisifes $H_f\geq (n-1+4k)c_0>(n-1+4k)c\geq 0$ for some constant $c_0>c\geq 0$. Let $d$ denote the distance function on $N$. Then
\begin{eqnarray}\label{dist5}
    \sup_{x\in N}d(x,\pt N)&\leq& \left\{\begin{array}{ll}
                                   \f {1}{c_0}, &\textrm{ if }c=0 \\
                                    & \\
                                  \f {1}c\coth^{-1}\frac{c_0}{c}, &\textrm{ if }c>0
                                \end{array}\right.
\end{eqnarray}
\end{prop}
\proof
As in the proof of diameter estimate \eqref{dist4}, we have
\begin{align}
    0\leq&\int_0^d\left((n-1)\varphi'(s)^2-\varphi(s)^2Ric_f(\gamma'(s),\gamma'(s))\right)ds\nonumber\\
    &\qquad-2\int_0^d\varphi(s)\varphi'(s)\frac d{ds}f(\gamma(s))ds-H_f(\gamma(d))\nonumber\\
    =&\int_0^d\left((n-1)\varphi'(s)^2-\varphi(s)^2Ric_f(\gamma'(s),\gamma'(s))\right)ds-H_f(\gamma(d))\nonumber\\
    &\qquad+2\int_0^d\frac d{ds}(\varphi(s)\varphi'(s))f(\gamma(s))ds-2\varphi(d)\varphi'(d)f(\gamma(d)),\label{eqA-1}
\end{align}
where in the second equality we used the integration by parts and the fact $\varphi(0)=0$.

If $c=0$, by choosing $\varphi(s)=s/d$ in \eqref{eqA-1}. Using the assumption $Ric_f\geq 0$ in $N$, $|f|\leq k$ and $H_f\geq (n-1+4k)c_0$ on $\pt N$, we have
\begin{eqnarray*}
  0 &\leq& \frac{n-1}d+\frac 2{d^2}\int_0^df(\gamma(s))ds-\frac 2df(\gamma(d))-(n-1+4k)c_0\\
  &\leq&\frac{n-1+4k}{d}-(n-1+4k)c_0.
\end{eqnarray*}
Therefore, we have
\begin{eqnarray}
  d&\leq &\frac{1}{c_0}.
\end{eqnarray}

If $c>0$, by choosing $\varphi(s)=\sinh(cs)/{\sinh(cd)}$ for $0\leq s\leq d$ in \eqref{eqA-1}. Using the assumption $Ric_f\geq -(n-1)c^2$ in $N$, $|f|\leq k$ and $H_f\geq (n-1+4k)c_0>(n-1+4k)c>0$ on $\pt N$, we have
\begin{eqnarray*}
  0 &\leq& (n-1)c\coth(cd)+\frac{2c^2}{\sinh^2(cd)}\int_0^d\cosh(2cs)f(\gamma(s))ds\\
  &&\qquad-2c\coth(cd)f(\gamma(d))-(n-1+4k)c_0 \\
  &\leq&(n-1+4k)c\coth(cd)-(n-1+4k)c_0,
\end{eqnarray*}
which is equivalent to
\begin{eqnarray}
  d &\leq& \frac 1c\coth^{-1}\frac{c_0}{c}
\end{eqnarray}
\endproof

\bibliographystyle{Plain}

\end{document}